\documentclass{article}
\usepackage{amsfonts}
\usepackage{amssymb}
\usepackage{latexsym}
\def\geqs{\geqslant}
\def\leqs{\leqslant}
\def\be{\begin{equation}}
\def\ee{\end{equation}}
\def\bea{\begin{eqnarray}}
\def\eea{\end{eqnarray}}
\def\uint{{\int_0^1}}
\def\uiint{{\uint\!\uint}}
\def\uidint{{\uint\!\cdots\!\uint}}
\def\idint{\mathop{\int\!\cdots\!\int}}
\def\pint{{\int_0^{\pi/2}}}
\def\pidint{{\pint\!\!\!\cdots\!\pint}}
\def\sumk{{\sum_{k=0}^\infty}}

\def\prez{\prec^{\phantom0}_0}
\def\sucz{\succ^{\phantom0}_0}
\newcommand{\R}{\mathbf{R}}
\newcommand{\Z}{\mathbf{Z}}
\newcommand{\Vol}{\mathop{\rm Vol}\nolimits}
\begin{document}
\centerline{\large On the Sums $\sum_{k=-\infty}^\infty (4k+1)^{-n}$}
\vspace*{2ex}
\centerline{Noam D.\ Elkies}

\vspace*{5ex}

\begin{quote}
{\bf Abstract.} The sum in the title is a rational multiple of~$\pi^n$
for all integers $n=2,3,4,\ldots$ for which the sum converges absolutely.
This is equivalent to a celebrated theorem of Euler.
Of the many proofs that have appeared since Euler,
a simple one was discovered only recently by Calabi~\cite{BKC}:
the sum is written as a definite integral over the unit \hbox{$n$-cube},
then transformed into the volume of a polytope $\Pi_n \subset \R^n$
whose vertices' coordinates are rational multiples of~$\pi$.
We review Calabi's proof, and give two further interpretations.
First we define a simple linear operator~$T$\/ on L$^2(0,\pi/2)$,
and show that $T$\/ is self-adjoint and compact,
and that $\Vol(\Pi_n)$ is the trace of $T^n$.
We find that the spectrum of~$T$\/ is $\{ 1/(4k+1) : k\in\Z \}$,
with each eigenvalue $1/(4k+1)$ occurring with multiplicity~$1$;
thus $\Vol(\Pi_n)$ is the sum of the \hbox{$n$-th} powers
of these eigenvalues.
We also interpret $\Vol(\Pi_n)$ combinatorially
in terms of the number of alternating permutations of $n+1$ letters,
and if $n$ is even also in terms
of the number of cyclically alternating permutations of $n$ letters.
We thus relate these numbers with $S(n)$
without the intervention of Bernoulli and Euler numbers
or their generating functions.
\end{quote}

{\bf 1.~INTRODUCTION.}
Euler proved that the sums
\be
\zeta(n) := \sum_{k=1}^\infty \frac1{k^n}
= 1 + \frac1{2^n} + \frac1{3^n} + \frac1{4^n} + \cdots
\label{zeta(n)}
\ee
for even $n\geqs2$ and
\be
L(n,\chi_4) := \sumk \frac{(-1)^k}{(2k+1)^n}
= 1 - \frac1{3^n} + \frac1{5^n} - \frac1{7^n} + - \cdots
\label{l4sum}
\ee
for odd $n\geqs3$ are rational multiples of~$\pi^n$.
(See Ayoub~\cite{Ayoub} for more
on Euler's work on these and related sums.)
This result can be stated equivalently as follows:
\be
S(n) := \sum_{k=-\infty}^\infty \frac1{(4k+1)^n}
\label{S(n)}
\ee
is a rational multiple of $\pi^n$ for all integers $n=2,3,4,\ldots$
for which the sum converges absolutely.
This is equivalent to~(\ref{zeta(n)}) and~(\ref{l4sum}) because
\be
S(n) = \cases{
(1 - 2^{-n}) \, \zeta(n),& if $n$ is even;\cr
L(n,\chi_4),& if $n$ is odd.
}
\label{S(n)=}
\ee
For future reference we tabulate for $n\leqs10$ the rational numbers
$\pi^{-n} S(n)$, as well as $\pi^{-n}\zeta(n)$ for $n$ even.
To this end, we define
\be
S(1) = 1-\frac13+\frac15-\frac17+-\cdots = \frac \pi 4.
\label{S(1)}
\ee
\begin{center}
\begin{tabular}{c|cccccc}
        $n$        &  1  &  2  &  3   &  4   &   5    &   6   \\ \hline
  $\pi^{-n} S(n)$  & 1/4 & 1/8 & 1/32 & 1/96 & 5/1536 & 1/960 \\ \hline
$\pi^{-n} \zeta(n)$&     & 1/6 &      & 1/90 &        & 1/945
\end{tabular}

\begin{tabular}{c|cccc}
        $n$        &     7     &      8    &      9      &     10 \\ \hline
  $\pi^{-n} S(n)$  & 61/184320 & 17/161280 & 277/8257536 & 31/2903040 \\ \hline
$\pi^{-n} \zeta(n)$&           &  1/9450   &             &  1/93555
\end{tabular}
\end{center}
One standard proof of the rationality of $\pi^{-n} S(n)$ is via
the generating function
\be
G(z) := \sum_{n=1}^\infty S(n) z^n
= \sum_{n=1}^\infty
  \left(\sum_{k=-\infty}^\infty \frac1{(4k+1)^n}\right)z^n,
\label{G(z)}
\ee
in which the inner sum is taken in the order
$k=0,-1,1,-2,2,-3,3,\ldots$.
The power series representing~$G(z)$ converges
for all $z$ such that $|z|<1$.
Since the sum of the terms with $n>1$ converges absolutely,
we may interchange the order of summation in~(\ref{G(z)}), obtaining
\be
\sum_{k=-\infty}^\infty
 \left( \sum_{n=1}^\infty \frac{z^n}{(4k+1)^n} \right)
= z \left(
 \frac1{1-z} - \frac1{3+z} + \frac1{5-z} - \frac1{7+z} + - \cdots
 \right).
\label{partfrac}
\ee
Comparing the latter sum with the partial-fraction expansions
for the tangent and cosecant, we find that
\be
G(z) = \frac{\pi z}{4}
\left( \sec\frac{\pi z}{2} + \tan\frac{\pi z}{2} \right).
\label{sec,tan}
\ee
Since the Taylor series of $z(\sec(z)+\tan(z))$ about $z=0$
has rational coefficients, it follows from (\ref{sec,tan})
that for each~$n$ the coefficient $S(n)$ of~$z^n$ in~$G(z)$
is a rational multiple of $\pi^n$.

This also lets us easily compute the rational numbers $\pi^{-n} S(n)$.
The numbers for even and odd $n$ come from the even and odd parts
$(\pi z/4) \tan(\pi z/2)$ and $(\pi z/4) \sec(\pi z/2)$ of $G(z)$.
In the literature these are usually treated separately,
and their coefficients are expressed in terms of 
Bernoulli and Euler numbers $B_n$ and $E_{n-1}$, respectively,
which are related to $S(n)$ by the formulas
\be
B_{2m} = \frac{(-1)^{m-1}}{2^{2m-1}\pi^{2m}} \zeta(2m)
       = (-1)^{m-1} \frac2{2^{2m}-1} \pi^{-2m} S(2m),
\label{B}
\ee
\be
E_{2m} = \frac{(-1)^m}{(2m)!\, 2^{2m+2}} \pi^{-(2m+1)} S(2m+1).
\label{E}
\ee
(The formulas (\ref{B}) and (\ref{E}) specify only
the Bernoulli and Euler numbers of positive even order.
The odd-order ones all vanish except for $B_1=-1/2$.
This definition, as well as the initial value $B_0=1$,
is needed for another important use of these numbers;
namely, the formula
\be
\sum_{k=1}^N k^{n-1} = \frac1n \sum_{m=1}^n {n \choose m} B_m N^{n-m}
\label{Bsum}
\ee
for summing powers of the first $N$\/ integers.)
A short table of Euler and Bernoulli numbers of even order follows:
\begin{center}
\begin{tabular}{c|c|c|c|c|c|c|c|c|c|c}
$B_0$ & $E_0$ & $B_2$ & $E_2$ & $B_4$ & $E_4$ & $B_6$ & $E_6$ &
 $B_8$ & $E_8$ & $B_{10}$ \\ \hline
  1   &   1   &  1/6  & $-1$  &$-1/30$&   5   & 1/42  & $-61$ &
$-1/30$ & 1385 & 5/66
\end{tabular}
\end{center}
It is known that all the Euler numbers are integers; we shall give
a combinatorial interpretation of their absolute values
$(-1)^m E_{2m}$ at the end of this paper.

The sums $S(n)$ continue to attract considerable interest
in mathematical disciplines
ranging from Fourier analysis to number theory.
Euler's formulas predate the year 1750,
and over the years since Euler's time,
many new proofs of the rationality of $S(n)/\pi^n$ have been given.
But it was only recently that Calabi\footnote{
  The only paper by Calabi that contains this proof is one co-authored
  with Beukers and Kolk \cite{BKC}.  Nevertheless, the proof is due
  to Calabi alone; Beukers and Kolk's contribution to~\cite{BKC}
  concerns other aspects of that paper.  They were shown this proof
  by Don Zagier, who also introduced me to it.
  Note that Kalman~\cite{Kalman} also writes that he first
  learned of this proof in a talk by Zagier.}
found a proof using only
the formula for change of variables of multiple integrals.
For instance, to prove that $\zeta(2)=\pi^2/6$, or equivalently
that $S(n)=\pi^2/8$ for $n=2$, Calabi argues as follows:
Write each term $(2k+1)^{-2}$ of the infinite sum in (\ref{zeta(n)})
as $\uiint (xy)^{2k} dx\,dy$, and thus rewrite that sum as
\bea
\sumk \frac1{(2k+1)^2}
&=& \sumk\, \uiint (xy)^{2k} dx\,dy \\
&=& \uiint \left( \sumk (xy)^{2k} \right) dx\,dy \label{eq0} \\
&=& \uiint \frac {dx\,dy} {1-(xy)^2} \label{squint}.
\eea
(The interchange of sum and integral in (\ref{eq0}) is readily
justified, for instance by observing the positivity of each
integrand.)  Then perform the change of variable
\be
x=\frac{\sin u}{\cos v},\qquad
y=\frac{\sin v}{\cos u},\label{cv2}
\ee
under which the integrand in (\ref{squint})
miraculously transforms to $1\,du\,dv$,
and the region of integration $\{ (x,y) \in \R^2 : 0<x,y<1\}$
is the one-to-one image of the isosceles right triangle
$\{ (u,v) \in \R^2: u,v>0,\;u+v<\pi/2\}$
(these assertions will be proved in greater generality below).
Thus the value of the integral (\ref{squint})
is just the area $\pi^2/8$ of that triangle, Q.E.D.

In general, Calabi writes $S(n)$ as a definite integral over
the \hbox{$n$-cube} $(0,1)^n$ and transforms it to the integral
representing the volume of the \hbox{$n$-dimensional} polytope
\be
\Pi_n := \{(u_1,u_2,\ldots,u_n)\in{\bf R}^n:
u_i>0,\; u_i+u_{i+1}<\frac\pi2\; (1\leqs i\leqs n)\}.
\label{polytope}
\ee
Note that the $u_i$ are indexed cyclically mod~$n$, so
\be
u_{n+1} := u_1,
\label{mod_n}
\ee
here and henceforth.
Since all the coordinates of each vertex of~$\Pi_n$
are rational multiples of~$\pi$,
the volume of~$\Pi_n$ must be a rational multiple of $\pi^n$.
It turns out that there is another way to interpret $S(n)$
as the volume of $\Pi_n$; this alternative approach requires
more analytical machinery, but better explains
the appearance of the sum $S(n)$.
We shall also give combinatorial interpretations of $S(n)$
by relating this volume to the enumeration of alternating permutations
of $n+1$ letters, and to the enumeration of cyclically alternating
permutations of $n$ letters when $n$ is even.
This leads to known formulas involving $B_n$ and $E_{n-1}$.
Our treatment via $\Pi_n$ and another polytope
relates those permutation counts directly to $S(n)$
without the intervention of Bernoulli and Euler numbers
or their generating functions.

To keep this paper self-contained, we first review
Calabi's transformation~\cite{BKC} that proves $S(n)=\Vol(\Pi_n)$.
This elegant proof remains little-known
(except possibly for the case $n=2$ shown above,
which was the second of Kalman's six proofs of $\zeta(2)=\pi^2/6$
\cite{Kalman}), and deserves wider exposure.
We then give the analytic interpretation of both $S(n)$ and
$\Vol(\Pi_n)$ as the trace of $T^n$ for a certain compact
self-adjoint operator~$T$\/ on the Hilbert space L$^2(0,\pi/2)$.
(I thank the referee for bringing to my attention a similar evaluation
of an integral studied by Kubilius~\cite{Kubilius}.)
Finally, we relate $S(n)$ and polytope volumes to the enumeration
of alternating and cyclically alternating permutations.

\vspace*{2ex}

{\bf 2.~EVALUATING $S(n)$ BY CHANGE OF VARIABLES.}
Following Calabi~\cite{BKC},
we generalize (\ref{cv2}) to the $n$-variable transformation
\be
x_1=\frac{\sin u_1}{\cos u_2},\quad
x_2=\frac{\sin u_2}{\cos u_3},\quad\ldots,\quad
x_{n-1}=\frac{\sin u_{n-1}}{\cos u_n},\quad
x_n=\frac{\sin u_n}{\cos u_1}, \label{cvn}
\ee
some of whose properties are established in the following two lemmas.

{\bf Lemma 1.}
{\em The Jacobian determinant of the transformation (\ref{cvn}) is
\be
\frac{\partial(x_1,\ldots,x_n)}{\partial(u_1,\ldots,u_n)}
= 1 \pm (x_1 x_2 \cdots x_n)^2,  \label{jacdet}
\ee
the sign $-$ or $+$ chosen according to whether $n$ is even or odd.}

{\em Proof}.
The partial derivatives $\partial x_i / \partial u_j$ are given by
\be
\frac{\partial x_i}{\partial u_j} = \cases{
 (\cos u_i) / (\cos u_{i+1}),& if $j=i$;\cr
 (\sin u_i \sin u_{i+1}) / (\cos^2 u_{i+1}),&
   if $j\equiv i+1\bmod n$;\cr
 0,& otherwise.\cr
 } \label{pderivs}
\ee
Thus the expansion of the Jacobian determinant
has only two nonzero terms,
one coming from the principal diagonal $j=i$,
one from the cyclic off-diagonal $j\equiv i+1\bmod n$.
The product of the principal diagonal entries
$(\cos u_i)/(\cos u_{i+1})$ simplifies to~1,
and always appears with coefficient $+1$.
The product of the off-diagonal entries is
\be
\prod_{i=1}^n \frac{\sin u_i \sin u_{i+1}}{\cos^2 u_{i+1}}
= \prod_{i=1}^n x_i^2,
\ee
and appears with coefficient $(-1)^{n-1}$, the sign of an
$n$-cycle in the permutation group.  Therefore the Jacobian
determinant is given by (\ref{jacdet}), as claimed.  $\Box$

{\bf Lemma 2.}
{\em The transformation (\ref{cvn}) maps the polytope $\Pi_n$
one-to-one to the open unit cube $(0,1)^n$.}

{\em Proof}.
Certainly if $u_i,u_{i+1}$ are positive and $u_i+u_{i+1}<\pi/2$ then
\be
0 < x_i = \frac{\sin u_i}{\cos u_{i+1}}
< \frac{\sin(\pi/2-u_{i+1})}{\cos u_{i+1}} = 1. \label{ineq}
\ee
Likewise we see that, given arbitrary $x_i$ in~$(0,1)$,
any $(u_1,\ldots,u_n)$ in~$(0,\pi/2)^n$ satisfying (\ref{cvn})
must lie in $\Pi_n$.  It remains to show that there exists
a unique such solution $(u_1,\ldots,u_n)$.
Rewrite the equations (\ref{cvn}) as
\be
u_1=f_{x_1}(u_2),\
u_2=f_{x_2}(u_3),\ \ldots,\
u_{n-1}=f_{x_{n-1}}(u_n),\
u_n=f_{x_n}(u_1), \label{nvc}
\ee
where $f_x$ ($0<x<1$) is the map
\be
f_x(u):=\sin^{-1}(x\cos u) \label{fxdef}
\ee
of the interval $(0,\pi/2)$ to itself.  Since
\be
\left| \frac{d}{du} f_x(u) \right| =
\left| -\frac{x\sin u}{\sqrt{1-x^2\cos^2 u}} \right| <
\frac{x\sin u}{\sqrt{x^2-x^2\cos^2 u}} = 1,
\label{contr}
\ee
each $f_{x_i}$ is a contraction map of
the interval $[0,\pi/2]$; hence so is their composite, which
thus has a unique fixed point in that interval.
This point cannot be at either endpoint, because $f_x(\pi/2)=0$
and $f_x(0)=\sin^{-1}(x)$ belongs to $(0,\pi/2)$
for all $x$ in the interval $(0,1)$.
Therefore
$f_{x_1} \circ f_{x_2} \circ\,\cdots\,\circ f_{x_{n-1}} \circ f_{x_n}$
has a unique fixed point in $(0,\pi/2)$, and
the simultaneous equations (\ref{nvc}) have a unique solution
$(u_1,\ldots,u_n)$ for each $(x_1,\ldots,x_n)$, as claimed.  $\Box$

Thus we see that the volume of the polytope $\Pi_n$ is
\bea
\!\!\!\!\!\!\!\!
\idint_{u_i>0\atop u_i+u_{i+1}<\pi/2} \!1\;du_1\cdots du_n
& \! = \! & \uidint \frac{dx_1\cdots dx_n}{1\pm(x_1\cdots x_n)^2}
\label{eq1} \\
& \! = \! &
\uidint \sumk (-1)^{nk} (x_1\cdots x_n)^{2k}\, dx_1\cdots dx_n.
\label{eq2}
\eea
Note that when $n$ is even, the second integral in (\ref{eq1})
is improper due to the singularity at $(x_1,\ldots,x_n)=(1,\ldots,1)$,
but the change of variable remains valid because the integrand is
everywhere positive.  By absolute convergence we may now interchange
the sum and multiple integral in (\ref{eq2}), obtaining
\be
\sumk (-1)^{nk} \uidint (x_1\cdots x_n)^{2k}\, dx_1\cdots dx_n
=\sumk \frac{(-1)^{nk}}{(2k+1)^n} = S(n).
\ee

We have thus proved:

{\bf Theorem 1.}
{\em The volume of the polytope $\Pi_n$ is $S(n)$ for all $n\geqs2$.}

{\bf Corollary 1.1.}
{\em $S(n)$ is a rational multiple of $\pi^n$ for all $n\geqs2$.}

Indeed, the volume of $\Pi_n$ is $(\pi/2)^n$ times the volume
of the polytope
\be
\frac2{\pi} \Pi_n = \{(v_1,v_2,\ldots,v_n)\in{\bf R}^n:
v_i>0,\; v_i+v_{i+1}<1\; (1\leqs i\leqs n)\},
\label{rat.polytope}
\ee
which is clearly a rational number.  $\Box$

\textbf{Remark}.  These results hold also when $n=1$,
though a bit more justification is needed because the alternating sum
(\ref{l4sum}) no longer converges absolutely when $n=1$.
In that case $\Pi_n$ reduces to the line segment $0<u_1<\pi/4$,
and the change of variable (\ref{cvn}) simplifies to $x_1=\tan u_1$,
so we recover the evaluation of
\be
S(1)=1-\frac13+\frac15-\frac17+-\cdots
\ee
as the arctangent integral
\be
\uint\frac{dx}{1+x^2}=\tan^{-1}(1)=\frac\pi4.
\ee

\vspace*{2ex}

{\bf 3.~RELATING $S(n)$ TO $\Pi_n$ VIA LINEAR OPERATORS.}
For $u,v$ in $(0,\pi/2)$, define $K_1(u,v)$ to be
the characteristic function of the isosceles right triangle
$\{ (u,v) \in \R^2: u,v>0,\;u+v<\pi/2\}$
encountered in the introduction, that is,
\be
K_1(u,v)=\cases{1,& if $u+v<\pi/2$;\cr 0,& otherwise.\cr}
\label{k1def}
\ee
We may then rewrite the volume of the polytope $\Pi_n$ as
\bea
&&
\pidint \prod_{i=1}^n K_1(u_i,u_{i+1})
\,du_1\,du_2\cdots du_n \label{cubint}\\
&&=\pint K_n(u,u)\, du
\label{trint}
\eea
(recall that $u_{n+1} = u_1$), where
\bea
&&\!\!\!\!\!\!\!\!\!\!
K_n(u,v)=\pint K_1(u,u_1) K_{n-1}(u_1,v)\, du_1 \label{kinduct}\\
\!\!\!\!\!\!\!\!=&&\!\!\!\!\!\!\!\!\!\!
\pidint
K_1(u,u_1) \cdot \prod_{i=1}^{n-2}K_1(u_i,u_{i+1}) \cdot K_1(u_{n-1},v)
\,du_1\,du_2\cdots du_{n-1} \ \label{kint}
\eea
(the equivalence of these two formulas, and thus also of (\ref{cubint})
with (\ref{trint}), is easily established by induction on~$n$).
We now interpret $K_n$ and the integral (\ref{trint}) in terms
of linear operators on L$^2(0,\pi/2)$.
Let $T$\/ be the linear operator with kernel $K_1(\cdot,\cdot)$
on L$^2(0,\pi/2)$:
\be
(Tf)(v) = \pint f(u) K_1(u,v)\, du = \int_0^{(\pi/2)-v} f(u)\, du.
\label{tdef}
\ee
Then we see from either (\ref{kinduct}) or (\ref{kint}) that
$K_n(\cdot,\cdot)$ is the kernel of $T^n$:
\be
(T^n\!f)(v) = \pint f(u) K_n(u,v)\, du.
\label{tnint}
\ee

The next lemma gives the spectral decomposition of this operator $T$,
and thus also of its powers $T^n$.

{\bf Lemma 3.}
{\em The transformation $T$\/ is a compact, self-adjoint operator
on {\rm L\/}$^2(0,\pi/2)$.  Its eigenvalues, each of multiplicity one,
are $1/(4k+1)\;\; (k\in\bf Z)$; the \hbox{corresponding} orthogonal
eigenfunctions are $\cos ((4k+1)u)$.}

{\em Proof}.
$T$\/ is self-adjoint because its kernel is symmetric:
$K_1(u,v)=K_1(v,u)$.  Compactness can be obtained from general
principles (the functions in the image $\{Tf:\|f\|\leqs1\}$
of the unit ball are a uniformly continuous family), or from the
determination of $T$\/'s eigenvalues.  Let $\lambda$, then, be
an eigenvalue of $T$, and $f$\/ a corresponding eigenfunction, so
\be
\int_0^{(\pi/2)-v} f(u)\,du = \lambda f(v) \label{eigen}
\ee
for almost all $v$ in $(0,\pi/2)$.  Note that $\lambda$ may not vanish,
because then (\ref{eigen}) would give $f=0$ in L$^2(0,\pi/2)$.
So we may divide (\ref{eigen})
by~$\lambda$, and use the left-hand side to realize $f$\/ as a
continuous function, and again to show that it is differentiable,
with
\be
f(\frac\pi2-v) = -\lambda f'(v) \label{dfq1}
\ee
for all $v$ in $(0,\pi/2)$.
Differentiating (\ref{dfq1}) once more, we find that
\be
\lambda^2 f''(v) = -\lambda \frac{d}{dv} f(\frac\pi2-v) =
\lambda f'(\frac\pi2-v) = -f(v), \label{dfq2}
\ee
whence
\be
f(v)=A\cos\frac{v}{\lambda} + B\sin\frac{v}{\lambda} \label{dfqsol}
\ee
for some constants $A$ and~$B$.  But from (\ref{eigen}) we see that
$f(\pi/2)=0$; substituting this into (\ref{dfq1}) we obtain
$f'(0)=0$, so $B=0$.  The condition $f(\pi/2)=0$ then becomes
$\cos(\pi/2\lambda)=0$ and forces $\lambda$ to be the reciprocal
of an an odd integer, say $\lambda=1/m$.  Now $\lambda$ and $-\lambda$
would both give rise to the same function $f(v)=\cos(v/\lambda)$,
but only one of them may be its eigenvalue.  To choose the sign,
take $v=0$ in (\ref{eigen}), finding that
\be
\lambda=\lambda f(0)=\pint\cos(mu)\,du=\frac1m\sin\frac{m\pi}{2}
=\frac{(-1)^{(m-1)/2}}{m},  \label{sign}
\ee
or equivalently that $\lambda=1/m$ with $m\equiv1\bmod4$.
We then easily confirm that each of these
$\lambda=1,-1/3,1/5,-1/7,\ldots$ and the corresponding
$f(v)=\cos(v/\lambda)$ satisfy (\ref{eigen}) for all
$v$ in $(0,\pi/2)$, completing the proof of the Lemma.
Alternatively, having obtained the eigenfunctions
$\cos(u)$, $\cos(3u)$, $\cos(5u),\ldots$, we need only invoke the
theory of Fourier series to show that these form an orthogonal
basis for L$^2(0,\pi/2)$ and then verify that they satisfy (\ref{eigen})
with the appropriate $\lambda$.  $\Box$

{\bf Corollary 3.1.}
{\em
The transformation $T^n$ is a compact self-adjoint operator on
L$^2(0,\pi/2)$.  Its eigenvalues, each of multiplicity one, are
$1/(4k+1)^n$ ($k\in\bf Z$), with corresponding
eigenfunctions $\cos((4k+1)u)$.
}

In particular, once $n\geqs2$, the sum of the eigenvalues of $T^n$
converges absolutely, so $T^n$ is of {\sl trace class}
(see Dunford and Schwartz \cite[XI.8.49, pp.1086--7]{DS}),
and its {\sl trace} is the sum
\be
\sum_{k=-\infty}^\infty \frac1{(4k+1)^n} = S(n)
\label{trace}
\ee
of these eigenvalues.  But it is known that the trace of a
trace-class operator is the integral of its kernel over the diagonal
(a continuous analog of the fact that the trace of a matrix is the sum
of its diagonal entries \cite[XI.8.49(c),\hbox{pp.1086--7}]{DS}).
Thus the trace $S(n)$ of $T^n$ is given by the integral
(\ref{trint}), i.e., by the volume of~$\Pi_n$.
So we have an alternative proof of Theorem 1, in which the power sum
(\ref{trace}) appears naturally, without separating the cases
of even and odd~$n$.

For future use, we give the orthogonal expansion of an arbitrary
L$^2$ function and a consequence of Corollary~3.1:

{\bf Corollary 3.2.}
{\em
For any $f$\/ in~${\rm L}^2(0,\pi/2)$ we have
\be
f = \sum_{k=-\infty}^\infty f_k \cos((4k+1)u),
\label{sumfk}
\ee
with coefficients $f_k$ given by
\be
\frac4\pi \pint f(u) \cos((4k+1)u) \, du.
\label{fk}
\ee
For each $n\geqs0$ we have
\be
\pint f(u) \, (T^n f)(u) \, du
= \frac\pi4 \sum_{k=-\infty}^\infty \frac{f_k^2}{(4k+1)^n}.
\label{f,Tf}
\ee
}

{\em Proof}.
Formulas (\ref{sumfk}) and (\ref{fk}) follow as usual from
the orthogonality of the eigenfunctions $\cos((4k+1)u)$
and the fact that
\be
\pint \cos^2((4k+1)u) \, du = \frac\pi 4
\label{pi/4}
\ee
for each integer~$k$.  This, together with the eigenvalues of~$T^n$
given in Corollary~3.1, yields (\ref{f,Tf}) as well.~~$\Box$

\vspace*{2ex}

{\bf 4.~ALTERNATING PERMUTATIONS, $S(n)$, AND $\Pi_n$.}
We shall see that $S(n)$ is closely related with the enumeration
of alternating (also known as ``up-down'' or ``zig-zag'') permutations
of $n$ letters.  A permutation $\sigma$ of
\be
[n]:=\{1,2,\ldots,n\}
\label{[n]}
\ee
is said to be {\sl alternating} if
\be
\sigma(1)<\sigma(2)>\sigma(3)<\sigma(4)>\,<\cdots
{{\lower1.5ex\hbox{$>$}}\atop{\raise1.5ex\hbox{$<$}}}\sigma(n)
\label{altdef}
\ee
(the last inequality is $<$ or $>$ for $n$ even and odd, respectively);
such $\sigma$ is {\sl cyclically alternating} if, in addition,
$n$ is even and $\sigma(n)>\sigma(1)$ \cite{Stanley}.\footnote{
  The term ``permutation altern\'ee'' was introduced by Andr\'e
  \cite{Andre} (see also \cite{Andre1})
  in the paper that first introduced alternating permutations
  and related their enumeration to the power series
  for $\sec x$ and $\tan x$.  This notion should not be confused with
  the ``alternating group'' of even permutations of~$[n]$.
  }
Let $A(n)$ be the number of alternating permutations of~$[n]$,
and $A_0(n)$ the number of cyclically alternating permutations
when $n$ is even.  We tabulate these numbers for $n \leqs 10$:
\centerline{
\begin{tabular}{c|cccccccccc}
  $n$   & 1 & 2 & 3 & 4 &  5 &  6 &  7  &   8  &   9  &  10 \\ \hline
$A(n)$  & 1 & 1 & 2 & 5 & 16 & 61 & 272 & 1385 & 7936 & 50521 \\ \hline
$A_0(n)$&   & 1 &   & 4 &    & 48 &     & 1088 &      & 39680
\end{tabular}
}
The table suggests a relationship between $A(2m-1)$ and $A_0(2m)$,
which we prove next.

{\bf Lemma 4.} {\em For all $m=1,2,3,\ldots$, the counts
$A(2m-1)$ and $A_0(2m)$ are related by
\be
A_0(2m) = m A(2m-1).
\label{A0,A}
\ee
}%

{\em Proof}.
We construct a \hbox{1-to-$m$} correspondence
between alternating permutations of~$[2m-1]$
and cyclically alternating permutations of~$[2m]$.
Note first that,
if $\sigma$ is a cyclically alternating permutation of~$[2m]$,
then so is
\be
\bigl(\sigma(2j+1),\sigma(2j+2),\ldots,\sigma(2j+2m-1)\bigr)
= \sigma \circ \tau^{2j}
\label{sigmacyc}
\ee
for each $j=0,1,2,\ldots,m-1$, where $\tau$ is the \hbox{$2m$-cycle}
sending each $i$ to $i+1$ (and as usual $i+1$ and $2j+1,2j+2,\ldots$
are taken mod~$2m$).  This partitions the cyclically alternating
permutations into sets of~$m$.   Now each of these sets contains
a unique permutation~$\sigma_0$ taking $2m$ to~$2m$.
But such $\sigma_0$ correspond bijectively to
the alternating permutations
\hbox{
$\bigl( \sigma_0(1), (\sigma_0(2), \ldots, \sigma_0(2m-1) \bigr)$%
}
of~$[2m-1]$.
This establishes the identity (\ref{A0,A}).~~$\Box$

The number $A_0(n)$ of cyclically alternating permutations of $[n]$
is given by the formula
\be
A_0(n) = 2^{n-1}(2^n-1)|B_n|
\label{altenum}
\ee
(see, for instance, Stanley \cite{Stanley}).
By (\ref{B}), this is equivalent to
\be
A_0(n) = \left(\frac2\pi\right)^{\!n}\! n!\, S(n).
\label{altenum_B}
\ee
The identity (\ref{altenum}) is usually obtained by identifying
the exponential generating functions, not via $S(n)$ and the zeta
function.  But $A_0(n)$ can also be expressed directly
in terms of the volume of the polytope $\Pi_n$ (or rather
in terms of the scaled polytope (\ref{rat.polytope})),
thus leading to (\ref{altenum_B}) and showing
that the appearance of $S(n)$ there is not merely accidental.

A general principle for enumerating permutations \cite{Stanley}
shows that the number of cyclically alternating permutations of $[n]$
is $n!$ times the volume in the unit cube $(0,1)^n$
of the region $P_0(n)$ determined by the inequalities
\be
t_1<t_2>t_3<t_4>\cdots<t_n>t_1.
\label{altreg}
\ee
(Recall that ``cyclically alternating permutations of $[n]$''
can exist only if n is even.)
This is because $P_0(n)$ is the \textit{order polytope} associated
to the partial order $\prez$ on~$[n]$ in which
\be
1 \prez 2 \sucz 3 \prez 4 \sucz \cdots \prez n \sucz 1
\label{prez}
\ee
and all other pairs in $[n]$ are incomparable.  In general,
the (open) order polytope
associated by Stanley to a partial order $\prec$ on~$[n]$
is the set of all $(t_1,\ldots,t_n)$ in the unit cube such that
$t_i<t_j$ whenever $i\prec j$; and the volume of this polytope
is an integer multiple of~$1/n!$.  Specifically, we cite
Corollary~4.2 from~\cite{Stanley}.

{\bf Lemma 5.} {\em The volume of the order polytope
associated to any partial order $\prec$ on~$[n]$
is $1/n!$ times the number of permutations~$\sigma$ of~$[n]$
such that $\sigma(i)<\sigma(j)$ whenever $i\prec j$.}

In other words, the volume is $1/n!$ times the number
of extensions of~$\prec$ to a linear order on~$[n]$.
To see that this is equivalent to the assertion of the lemma,
consider that there are $n!$ linear orders on~$[n]$,
each determined by the permutation of~$[n]$
that sends the the minimal element to~$1$, the next one to~$2$,
and so on until the maximal element is sent to~$n$.
This order extends~$\prec$
if and only if $\sigma(i)<\sigma(j)$ whenever $i\prec j$.

Lemma 5 appears in~\cite{Stanley}
as part of a corollary to a much more powerful theorem.
For our purposes the following direct proof suffices.

{\em Proof}.  Decompose the closed unit $n$-cube into $n!$ simplices,
one for each permutation~$\sigma$ of~$[n]$, so that the simplex
indexed by~$\sigma$ consists of all $t_1,\ldots,t_n$ in $[0,1]$ with
\be
t_{\sigma^{-1}(1)} \leqs t_{\sigma^{-1}(2)}
 \leqs \cdots \leqs t_{\sigma^{-1}(n)}.
\label{order}
\ee
Each of these simplices has the same volume;
hence this common volume is $1/n!$.
The union of those simplices indexed by $\sigma$
satisfying $\sigma(i)<\sigma(j)$ whenever $i\prec j$
is the closure of the order polytope of~$\prec$.
Thus the volume of this polytope is $1/n!$
times the number of such~$\sigma$.~~$\Box$

Equivalently, and perhaps more intuitively:
the volume of the order polytope is the probability that
$n$ independent variables~$t_i$ drawn uniformly at random from $[0,1]$
satisfy $t_i<t_j$ whenever $i\prec j$, i.e., that their linear order
inherited from~$[0,1]$ extends the partial order~$\prec$.
This, however, is the same as the probability
that a randomly chosen permutation of~$[n]$
yields a partial order extending~$\prec\,$,
because (excepting the negligible case that some $t_i$ coincide)
the order of the~$t_i$ determines a permutation, and all permutations
are equally likely.

We return now to the order polytope~(\ref{altreg})
associated to~$\prec_0$.  The affine change of variables
\be
u_i=\cases{t_i,& $i$ odd;\cr 1-t_i,& $i$ even,\cr} \label{switch}
\label{cov}
\ee
which has constant Jacobian determinant $(-1)^{n/2}$,
transforms this region (\ref{altreg}) to the familiar polytope
$v_i>0,\,v_i+v_{i+1}<1$,\footnote{
  Richard Stanley notes that this polytope~(\ref{rat.polytope})
  is also a special case of a construction from his
  paper~\cite{Stanley}: it is the \textit{chain polytope}
  associated with the same partial order~$\prez$.  The result
  in~\cite{Stanley} that contains our Lemma~4 asserts that
  the chain and order polytopes associated with any partial order
  have the same volume.  In general this is proved by a piecewise
  linear bijection, but for partial orders of ``rank~1''$\!$,
  i.e., for which there are no distinct $a,b,c$ such that
  $a \succ b \succ c$, the polytopes are equivalent by a single
  affine chain of variables.  The partial order $\prez$ has
  rank~1, as does the partial order we define next in~(\ref{prec})
  to deal with $A(n)$.  Stanley's affine change of variables 
  for these two partial orders is just our~(\ref{cov});
  thus this part of our argument is again a simple special case of his.
  }
whose volume we have already identified with $(2/\pi)^n S(n)$
in two different ways.  Thus
\be
A_0(n) = \left(\frac2\pi\right)^{\!n} \; n!\ S(n)
\label{A_0(n)}
\ee
for all even~$n$.

By Lemma~4, we also recover a formula for $A(2m-1)$:
\be
A(2m-1) = \frac{A_0(2m)}{m} = \frac{2^{2m-1}(2^{2m}-1)}{m} |B_{2m}|
        = \frac{2^{2m+1}}{\pi^{2m}} (2m-1)!\, S(2m).
\label{A(2m-1)}
\ee
In other words,
\be
A(n) = \frac{2^{n+2}}{\pi^{n+1}} \, n!\ S(n+1)
\label{A(n)}
\ee
when $n$ is odd.  We next prove this formula directly for all $n$,
whether even or odd.

{\bf Theorem 2.}
{\em
The number $A(n)$ of alternating permutations of~$[n]$
is given by (\ref{A(n)}) for every positive integer~$n$.
}

{\em Proof}.
Let $\prec$ be the partial order on $[n]$ in which
\be
1 \prec 2 \succ 3 \prec 4 \succ \prec \cdots
{{\lower2ex\hbox{$\succ$}}\atop{\raise2ex\hbox{$\prec$}}} n
\label{prec}
\ee
and all other pairs in $[n]$ are incomparable.
Then $A(n)$ is the number of permutations~$\sigma$ of~$[n]$
such that $\sigma(i)<\sigma(j)$ whenever $i\prec j$.  Accordingly,
$A(n)$ is $n!$ times the volume of the associated order polytope
\be
\{(t_1,t_2,\ldots,t_n)\in{\bf R}^n:
0<t_i<1,\; t_1<t_2>t_3<t_4>\,<
\cdots {{\lower1.2ex\hbox{$>$}}\atop{\raise1.2ex\hbox{$<$}}}t_n
\}.
\label{polytope2}
\ee
The change of variables (\ref{cov}) transforms (\ref{polytope2})
into the region
\be
\{(v_1,v_2,\ldots,v_n)\in{\bf R}^n:
v_i>0,\; v_i + v_{i+1} < 1 \; (1\leqs i\leqs n-1)\}.
\label{region2}
\ee
(N.B. This looks like the familiar $(2/\pi) \Pi_n$,
but in fact strictly contains $(2/\pi) \Pi_n$,
because we do not impose the condition $v_n+v_1<1$.)
On the other hand, under the further linear change of variable
$v_i = (2/\pi) u_i$, the region~(\ref{region2}) maps to
\be
\{(u_1,u_2,\ldots,u_n)\in{\bf R}^n:
u_i>0,\; u_i + u_{i+1} < \pi/2 \; (1\leqs i\leqs n-1)\},
\label{region2'}
\ee
a region whose volume is $(\pi/2)^n$ times larger
than that of~(\ref{region2}).  Thus the volume of (\ref{region2}) is
\bea
& & \left(\frac2\pi\right)^{\!n}
  \pidint \; \prod_{i=1}^{n-1} K_1(u_i,u_{i+1}) \, du_1 \cdots du_n
\nonumber \\
& \!\! = \!\! & \left(\frac2\pi\right)^{\!n}
  \pidint K_{n-1} (u_1,u_n) \, du_1 \cdots du_n.
\label{ipn-1}
\eea
Now by (\ref{tnint}) the function
\be
u_n \mapsto \pint K_{n-1} (u_1,u_n) \, du_1
\label{Tn-1(1)}
\ee
in ${\rm L}^2(0,\pi/2)$ is the image under $T^{n-1}$
of the constant function~$1$.
Thus the integral (\ref{ipn-1}) is $(2/\pi)^n$ times the inner product
of~$1$ and $T^{n-1}1$ in ${\rm L}^2(0,\pi/2)$.
By (\ref{f,Tf}) in Corollary~3.2, this inner product is
\be
\frac\pi4 \sum_{k=-\infty}^\infty \frac{c_k^2}{(4k+1)^{n-1}} \, ,
\label{1,Tn-1(1)}
\ee
where $c_k$ is the coefficient of $\cos((4k+1)u)$
in the orthogonal expansion of~$1$.
Using (\ref{fk}) of the same corollary, we calculate
\be
c_k = \frac4\pi \pint \cos((4k+1)u) \, du
= \frac4\pi \left.\frac{\sin((4k+1)u)}{4k+1}\right|_0^{\pi/2}
= \frac4\pi \, \frac1{4k+1}.
\label{ck}
\ee
Therefore,
\be
\frac{A_n}{n!} =
\left(\frac2\pi\right)^{\!n} \left(\frac4\pi\right) S(n+1),
\ee
from which (\ref{A(n)}) follows.~~$\Box$

\textbf{Remark}.
Some time before publishing \cite{Stanley}, Stanley proposed
the computation of the volume of the region in~(\ref{region2'})
as a \textsc{Monthly} problem~\cite{Prob2701}.  The published
solution~\cite{Sol2701} used generating functions to express
the volume in terms of $E_n$ or $B_{n+1}$.  But the
``Editor's Comments'' at the end of the solution include the note:
``Several solvers observed that [$n!$ times the volume of the polytope]
is the number of zig-zag permutations of $1,2,\ldots,n$, \ldots'' \
This suggests that, even if the elementary proof of Lemma~4
is not yet in the literature, it is obvious enough that
these solvers at least implicitly recognized it, and applied it
together with the change of variables~(\ref{cov})
to relate the polytope volume to~$A_n$.

Theorem~2, together with Lemma~4, yields the following amusing
corollary: as $m\rightarrow\infty$,
\be
\frac{A_0(2m)}{A(2m)} 
= m \frac{A(2m-1)}{A(2m)} 
= \frac\pi4 \frac{S(2m)}{S(2m+1)} 
\rightarrow \frac\pi4.
\label{lim}
\ee
That is, a randomly chosen alternating permutation~$\sigma$ of~$[2m]$
is cyclically alternating with probability approaching $\pi/4$\,!
The convergence is quite rapid, with error falling
as a multiple of $3^{-2m}$; for instance, for $m=5$ we already find
$39680/50521 \approx 0.785416$, while $\pi/4 \approx 0.785398$.

When $n$ is even, say $n=2m$, the formula (\ref{A(n)}) simplifies to
\be
A(2m) = (-1)^m E_{2m} = |E_{2m}|
\label{A(2m)}
\ee
by (\ref{E}).  We have thus given a combinatorial interpretation of
the positive integer $(-1)^m E_{2m}$, as promised in the Introduction: 
it is the number of alternating permutations of~$[2m]$.

{\bf ACKNOWLEDGEMENTS.}
Don Zagier first showed me Calabi's proof
and later provided the reference~\cite{BKC}.
Richard Stanley furnished references to and remarks on his work on
order polytopes and the enumeration of alternating permutations,
and later gave the reference~\cite{Andre}.
The referees drew my attention to the work of
Kubilius~\cite{Kubilius} and Kalman~\cite{Kalman}.
They made other suggestions that improved the paper.
I thank them for their respective contributions.
I am also grateful to the Mathematical Sciences Research Institute
for its hospitality during the completion of this paper,
and to the National Science Foundation and the Packard Foundation
for partial support during the preparation of this work.

\end{document}